\documentclass[12pt]{article}

\usepackage{amsmath,amssymb,amsthm}
\usepackage{amsfonts}
\usepackage[mathscr]{eucal}
\begin{document}
\newtheorem{thm}{Theorem}
\numberwithin{thm}{section}
\newtheorem{lemma}[thm]{Lemma}
\newtheorem{remark}{Remark}
\newtheorem{corr}[thm]{Corollary}
\newtheorem{proposition}{Proposition}
\newtheorem{theorem}{Theorem}[section]
\newtheorem{deff}[thm]{Definition}
\newtheorem{case}[thm]{Case}
\newtheorem{prop}[thm]{Proposition}
\numberwithin{equation}{section}
\numberwithin{remark}{section}
\numberwithin{proposition}{section}
\newtheorem{corollary}{Corollary}[section]
\newtheorem{others}{Theorem}
\newtheorem{conjecture}{Conjecture}\newtheorem{definition}{Definition}[section]
\newtheorem{cl}{Claim}
\newtheorem{cor}{Corollary}
\newcommand{\ds}{\displaystyle}

\newcommand{\stk}[2]{\stackrel{#1}{#2}}
\newcommand{\dwn}[1]{{\scriptstyle #1}\downarrow}
\newcommand{\upa}[1]{{\scriptstyle #1}\uparrow}
\newcommand{\nea}[1]{{\scriptstyle #1}\nearrow}
\newcommand{\sea}[1]{\searrow {\scriptstyle #1}}
\newcommand{\csti}[3]{(#1+1) (#2)^{1/ (#1+1)} (#1)^{- #1
 / (#1+1)} (#3)^{ #1 / (#1 +1)}}
\newcommand{\RR}[1]{\mathbb{#1}}
\thispagestyle{empty}
\begin{titlepage}
\title{\bf Isoperimetric-type inequalities for iterated Brownian motion in $\RR{R}^{n}$}
\author{Erkan Nane\thanks{ Supported in part by NSF Grant \# 9700585-DMS }\\
Department of Mathematics\\
Purdue University\\
West Lafayette, IN 47906 \\
enane@math.purdue.edu}
\maketitle
\begin{abstract}
\noindent {\it We extend generalized isoperimetric-type inequalities to iterated Brownian motion over several domains in $\RR{R}^{n}$. 
These kinds of inequalities imply in particular that for 
domains of finite volume, the exit distribution and moments of the first exit time for iterated Brownian motion are 
maximized with the ball $D^{*}$ centered at
the origin, which has the same volume as $D$.}

\end{abstract}

\end{titlepage}

\section{Introduction}

Isoperimetric inequalities for Brownian motion and several other processes are well-known results. These kinds of inequalities give information 
about the geometry of
the domain (see \cite{bandle}, \cite{banlame}, \cite{luttinger1}-\cite{luttinger3}, \cite{mendez}). Our aim in this paper is to extend 
these generalized isoperimetric inequalities to iterated Brownian motion (IBM) over several domains in
$\RR{R}^{n}$.

To define the iterated  Brownian motion $Z_{t}$
started at $z \in \RR{R}$,
 let $X_{t}^{+}$, $X_{t}^{-}$ and $Y_{t}$
 be three  independent
one-dimensional Brownian motions, all started at $0$. Two-sided
Brownian motion is defined by
\[ X_{t}=\left\{ \begin{array}{ll}
X_{t}^{+}, &t\geq 0\\
X_{(-t)}^{-}, &t<0.
\end{array}
\right. \] Then the iterated Brownian motion started at $z \in
\RR{R}$ is
\[ Z_{t}=z+X(Y_{t}),\ \ \    t\geq 0.\]

In $\RR{R}^{n}$, one requires $X^{\pm}$ to be independent
$n-$dimensional Brownian motions. This is the version of the
iterated Brownian motion due to Burdzy, see \cite{burdzy1}.

Properties of IBM and Brownian-time Brownian motion (BTBM), defined below, analogous to the
properties of Brownian motion have been studied extensively by
several authors (see \cite{burdzy1, burdzy2}, \cite{klewis,koslew}, \cite{nane, nane2}, \cite{xiao} and the references therein). 
One of the main differences between these iterated processes  and Brownian motion is
that they  are not Markov processes.  However, these
processes have connections
with the
 parabolic operator
$\frac{1}{8} \Delta ^{2} -\frac{\partial}{\partial t}$, as
described in \cite{allouba1,deblassie}.

We next give the background on isoperimetric inequalities that motivates 
us to study generalized isoperimetric inequalities for iterated Brownian motion.

Let $D\subset \RR{R}^{n}$ be a domain of finite volume, and denote
by $D^{*}$ the ball in $\RR{R}^{n}$ centered at the origin with
the same volume as $D$. The class of quantities related to the
Dirichlet Laplacian in $D$ which are maximized or minimized by the
corresponding quantities for $D^{*}$ are often called generalized
isoperimetric inequalities (C. Bandle \cite{bandle}).

 Probabilistically generalized isoperimetric inequalities
read as
\begin{equation}\label{isopheat2}
P_{z}[\tau_{D}>t]\ \leq \ P_{0}[\tau_{D^{*}}>t]
\end{equation}
for all $z\in D$ and all $t>0$, where $\tau_{D}$ is the first exit
time of the Brownian motion from the domain $D$ and $P_{z}$ is the
associated probability measure when this process starts at $z$.
These kinds of Isoperimetric inequalities imply eigenvalue
inequalities for Dirichlet Laplacian, Integrals of Green functions
of the Dirichlet Laplacian and the moments of the exit time of
Brownian motion (see \cite{luttinger1} - \cite{luttinger3}). These
inequalities  serve as motivation for our results in this paper.

Isoperimetric inequalities have been extended in many directions including Brownian motion over convex domains in $\RR{R}^{n}$
where, one fixes the inradius $R_{D}$ of $D$, instead of fixing
the volume. The inradius is the supremum of the radius
of all the balls contained in $D$. Indeed,  a special case of the
theorem proved  by R. Ba\~{n}uelos, R. Lata\l a,
and P. Me\`{n}dez in \cite{banlame}  is the following; suppose $D\subset \RR{R}^{2}$
is a convex domain of finite inradius $R_{D}$. Then
\begin{equation}\label{isopheat3}
P_{z}[\tau_{D}>t]\ \leq \ P_{0}[\tau_{S(D)}>t]\ =\
P_{0}[\tau_{I(D)}>t]
\end{equation}
for all $z\in D$ and all $t>0$ where, $I(D)=(-R_{D}, R_{D})$ and
$S(D)= \RR{R}\times I(D)$.

More recently, P. Me\`{n}dez \cite{mendez} proved a generalization
of inequality (\ref{isopheat3}) in $\RR{R}^{n}$. Suppose $D\subset
\RR{R}^{n}$ is a convex domain of finite inradius $R_{D}$. Then
\begin{equation}\label{isopheat4}
P_{z}[\tau_{D}>t]\ \leq \ P_{0}[\tau_{S(D)}>t]\ =\
P_{0}[\tau_{I(D)}>t]
\end{equation}
for all $z\in D$ and all $t>0$ where, $I(D)=(-R_{D}, R_{D})$ and
$S(D)= \RR{R}^{n-1}\times I(D)$.

P. Me\`{n}dez  proved in \cite[Theorem 5.2]{mendez} a sharper
version of inequality (\ref{isopheat3}) in $\RR{R}^{2}$. Suppose
$D\subset \RR{R}^{2}$ is a convex domain of finite inradius
$R_{D}$ and finite diameter $d_{D}$. Then
\begin{equation}\label{isopheat5}
P_{z}[\tau_{D}>t]\ \leq \ P_{0}[\tau_{C(D)}>t]
\end{equation}
for all $z\in D$ and all $t>0$ where, $C(D)=\left[ \RR{R}\times
I(D)\right]\cap B(0,d_{D}\ - \ R_{D})$.

In analogy with ordinary Brownian motion and diffusions, if $\tau
_{D}(Z)$ is the first exit time of IBM
from domain $D $, started at $z\in D$,
$$
(i.e. \ \ \tau _{D} (Z) =\inf \{t\geq 0: \  Z_{t} \notin D \}).
$$ Then
 $P_{z}[\tau _{D}(Z) > t]$ provides a measure of the lifetime of the process in $D$.

Let $D\subset \RR{R}^{n}$ be an open set and let $D_{*}$ be any of
the sets $D^{*}$, $S(D)$, $I(D)$ or $C(D)$. $D_{*}=D^{*}$ the open
ball of the same volume as $D$ when $D$ is of finite volume.
$D_{*}=S(D)=\RR{R}^{n-1}\times I(D)$ where $I(D)=(-R_{D}, R_{D})$
when $D$ is a convex domain of finite inradius $R_{D}$.
$D_{*}=C(D)=\left[ \RR{R}\times I(D)\right]\cap B(0,d_{D}\ - \
R_{D})$ when $D\subset \RR{R}^{2}$ is a convex set of finite
inradius $R_{D}$ and finite diameter $d_{D}$.

The following is the main result of this paper. It is an extension of generalized isoperimetric-type inequalities to 
iterated Brownian motion over several domains.

\begin{theorem}\label{gen-isop-burdzy-IBM}
Let $D\subset \RR{R}^{n}$ be an open set, of finite volume or a
convex set of finite inradius or a convex set of finite inradius
with finite diameter, then
\begin{equation}\label{isop-burdzy-IBM}
P_{z}[\tau_{D}(Z)>t]\ \leq \ P_{0}[\tau_{D_{*}}(Z)>t]
\end{equation}
\end{theorem}

Let $\phi$ be an increasing function. Multiplying the inequality
(\ref{isop-burdzy-IBM}) by the derivative of $\phi $ and integrating the
resulting inequality in time from zero to infinity, we obtain the
following corollary.
\begin{corollary}\label{increasing0}
Let $\phi$ be a nonnegative increasing function, then
\begin{equation}
E_{z}(\phi (\tau_{D}(Z)))\leq E_{0}(\phi
(\tau_{D_{*}}(Z)))
\end{equation}
\end{corollary}
 In particular Corollary \ref{increasing0} implies that for $p\geq 1$,

\begin{equation}
E_{z}( (\tau_{D}(Z)))^{p}\leq E_{0}(
(\tau_{D_{*}}(Z))^{p})
\end{equation}

Theorem \ref{gen-isop-burdzy-IBM}, and
its corollaries imply that the distribution functions, and the
moments of the exit times are maximized with $D^{*}$ when $D$ is
an open set of finite volume, with $S(D)$ when $D$ is a convex
domain of finite inradius and with $C(D)$ when $D$ is a convex
domain of finite inradius and finite diameter.

We also obtain a version of the above Theorem \ref{gen-isop-burdzy-IBM}
 for another closely
related process, the so called
 Brownian-time Brownian motion (BTBM).  To define this,
 let $X_{t}$ and $Y_{t}$
 be two  independent
one-dimensional Brownian motions, all started at $0$. BTBM is
defined to be $Z_{t}^{1}=x+X(|Y_{t}|)$. Properties of this process
and its connections to PDE's have been studied in
\cite{allouba1} and  \cite{koslew}.

In analogy with ordinary Brownian motion and diffusions, if $\tau
_{D}(Z^{1})$ is the first exit time of BTBM
from domain $D $, started at $z\in D$,
 $P_{z}[\tau _{D}(Z^{1}) > t]$ provides a measure of the lifetime of the process in $D$. Then we have

\begin{theorem}\label{gen-isop-IBM}
Let $D\subset \RR{R}^{n}$ be an open set, of finite volume or a
convex set of finite inradius or a convex set of finite inradius
with finite diameter, then
\begin{equation}\label{isopIBM11}
P_{z}[\tau_{D}(Z^{1})>t]\ \leq \ P_{0}[\tau_{D_{*}}(Z^{1})>t]
\end{equation}
\end{theorem}

Let $\phi$ be an increasing function. Multiplying the inequality
(\ref{isopIBM11}) by the derivative of $\phi $ and integrating the
resulting inequality  in time from zero to infinity, we obtain the
following corollary.
\begin{corollary}\label{increasing}
Let $\phi$ be a nonnegative increasing function, then
\begin{equation}
E_{z}(\phi (\tau_{D}(Z^{1})))\leq E_{0}(\phi
(\tau_{D_{*}}(Z^{1})))
\end{equation}
\end{corollary}
 In particular Corollary \ref{increasing} implies that for $p\geq 1$,

\begin{equation}
E_{z}( (\tau_{D}(Z^{1})))^{p}\leq E_{0}(
(\tau_{D_{*}}(Z^{1}))^{p})
\end{equation}

Theorem \ref{gen-isop-IBM}, and
its corollaries imply that the distribution functions, and the
moments of the exit times are maximized with $D^{*}$ when $D$ is
an open set of finite volume, with $S(D)$ when $D$ is a convex
domain of finite inradius and with $C(D)$ when $D$ is a convex
domain of finite inradius and finite diameter.

We use integration by parts and the double integral representation of
the exit distribution of the first exit time of IBM from domains in $\RR{R}^{n}$ to prove 
Theorems \ref{gen-isop-burdzy-IBM} and \ref{gen-isop-IBM} in the next section.

\section{Proofs of Main results}

If $D\subset \RR{R} ^{n}$  is an open set, write
\[
\tau_{D}^{\pm}(z)=\inf \{ t\geq 0:\ \ X_{t}^{\pm} +z \notin D\},\]
and if $I\subset \RR{R}$ is an open interval, write
\[
\eta _{I}=\eta (I)= \inf \{ t\geq 0 :\ \  Y_{t}\notin I\}.
\]
So we have in the notation of the introduction 
$$
P[\tau_{D}(z)>t]=P_{z}[\tau_{D}>t].
$$

Recall that $\tau _{D}(Z)$ stands for the first exit time of
iterated Brownian motion from  $D$.
As in DeBlassie
\cite[\S3.]{deblassie}, we have by the continuity of the paths for
$Z_{t}=z+X(Y_{t})$, if $f$ is the probability density of $\tau
_{D}^{\pm}(z)$
\begin{equation}\label{translation1}
P_{z}[\tau _{D}(Z) > t]=\int_{0}^{\infty} \!\int_{0}^{\infty}
 P_{0}[\eta _{(-u,v)}>t]  f(u) f(v)dvdu.
\end{equation}

To prove isoperimetric-type inequalities we need the following
Lemmas.
\begin{lemma}\label{increasing-isop}
For any $u,v\in ( 0,\infty) $

$$ 
\left( \frac{\partial }{\partial v}
 P_{0}[\eta _{(-u,v)}>t] \right)\geq 0,
$$
and
$$ 
\left( \frac{\partial }{\partial u}
 P_{0}[\eta _{(-u,v)}>t] \right)\geq 0
$$

\end{lemma}

\begin{proof}
Fix  $u\in (0,\infty)$, let $0< v_{1}\leq v_{2} <\infty$. We have $(-u, v_{1})\subset (-u,v_{2})$, hence
$\eta _{(-u,v_{1})}\leq  \eta _{(-u,v_{2})}$. So, for any $t>0$, $ P_{0}[\eta _{(-u,v_{1})}>t]\leq P_{0}[\eta _{(-u,v_{2})}>t]$. Hence the first inequality holds.
The last inequality is proved similarly.
 
\end{proof}

\begin{lemma}[Integration by parts lemma]\label{integration-by-parts-lemma}
Let $D$ be a domain such that $\lim_{t\to\infty}P_{z}[\tau_{D}>t]=0$.
Then

$P_{z}[\tau _{D}(Z) > t] $ is equal to
\begin{eqnarray}
&= &\int_{0}^{\infty} \!\int_{0}^{\infty}
 P_{0}[\eta _{(-u,v)}>t]  f(u) f(v)dvdu.\label{lemma-first}\\
 & = &\int_{0}^{\infty} \!\int_{0}^{\infty}  \left( \frac{\partial }{\partial v}
 P_{0}[\eta _{(-u,v)}>t] \right) P[\tau_{D}(z)  >v]f(u)dvdu\label{lemma-int-by-parts1}\\
 & = &  \int_{0}^{\infty} \!\int_{0}^{\infty}  \left( \frac{\partial }{\partial u}
 P_{0}[\eta _{(-u,v)}>t] \right) P[\tau_{D}(z)  >u]f(v)dudv. \label{lemma-int-by-parts2}
\end{eqnarray}

\end{lemma}

\begin{proof}
Equation (\ref{lemma-int-by-parts1}) follows by integration by parts and observing  
$$ -\frac{\partial}{\partial v}P[\tau_{D}(z)  >v]=f(v)$$ and
$$
\lim_{u\rightarrow 0}P_{0}[\eta _{(-u,v)} > t] P[\tau_{D}(z) > u]=\lim_{u\rightarrow \infty }P_{0}[\eta _{(-u,v)} > t] P[\tau_{D} (z)> u]=0.
$$
Equation (\ref{lemma-int-by-parts2}) follows by Fubini-Tonelli theorem in equation (\ref{lemma-first}) and
 by similar arguments as in the proof of equation (\ref{lemma-int-by-parts1}).
\end{proof}
If we use integration by parts one more time, we obtain

$P_{z}[\tau _{D}(Z) > t]$
\begin{equation}\label{eqddo}
=\int_{0}^{\infty} \!\int_{0}^{\infty}  \left( \frac{\partial }{\partial u} \frac{\partial }{\partial v}
 P_{0}[\eta _{(-u,v)}>t] \right) P[\tau_{D}(z) >u]P[\tau_{D}(z)  >v]dvdu.
\end{equation}

We now can try to use the isoperimetric-type inequalities for Brownian motion to prove Theorem \ref{gen-isop-burdzy-IBM}. This is valid if one can show that
$$
 \frac{\partial }{\partial u} \frac{\partial }{\partial v}
 P_{0}[\eta _{(-u,v)}>t]\geq 0
$$
for all $u,v,t\in (0,\infty)$.
But this is not easy to show (maybe not true), since this is given by an infinite series given in \cite{nane}:
\begin{eqnarray}
\ &\ & \frac{\partial }{\partial u} \frac{\partial }{\partial v} P_{0}[\eta _{(-u,v)} >t]\nonumber\\
\ &\ & = 4 \ \sum_{n=0}^{\infty} \  \exp (-\frac{(2n+1)^{2}\pi ^{2}}{2(u+v)^{2}} t)
\left\{ \left[  \ \sin \ \frac{(2n+1)\pi u}{u+v}\ (\frac{1}{(u+v)^{4}}) \right. \right. \nonumber \\
 \ &\ & \left. \left.\ \ \ \ \ \ \  \times \ (\frac{\pi ^{3}(2n+1)^{3} t^{2}}{(u+v)^{2}} - \ 3\pi (2n+1)t \ + \ (2n+1)\pi uv)\ \
\right] \right.\nonumber\\
 \ &\ & +  \left. \left[ \ \cos\ \frac{(2n+1)\pi u}{u+v}\ (\frac{1}{(u+v)^{3}})(\frac{(2n+1)^{2}\pi ^{2}t (v-u)}{(u+v)^{2}} +\ \ u\ -\ v) \right] \right\}.
  \nonumber
\end{eqnarray}

So we follow another path to prove isoperimetric inequalities.

\begin{proof}[Proof of Theorem \ref{gen-isop-burdzy-IBM}]
The idea of the proof is using integration by parts and the corresponding generalized isoperimetric inequalities for Brownian motion.
Let $f^{*}$ denote the probability density of $\tau_{D_{*}}$.
By equation (\ref{translation1}) and integration by parts
 $ P_{z}[\tau_{D}(Z)>t]$ equals
\begin{eqnarray}
 &  &  \int_{0}^{\infty} \!\int_{0}^{\infty}  \left( \frac{\partial }{\partial v}
 P_{0}[\eta _{(-u,v)}>t] \right) P[\tau_{D}(z)  >v]f(u)dvdu\label{int-by-parts0}\\
 &\leq & \int_{0}^{\infty} \!\int_{0}^{\infty}  \left( \frac{\partial }{\partial v}
 P_{0}[\eta _{(-u,v)}>t] \right) P[\tau_{D_{*}}(0)  >v]f(u)dvdu \label{first-bm-isop}\\
 & = &\int_{0}^{\infty} \!\int_{0}^{\infty}
 P_{0}[\eta _{(-u,v)}>t]  f(u) f^{*}(v)dvdu\label{int-by-parts2}\\
 & = & \int_{0}^{\infty} \!\int_{0}^{\infty}  \left( \frac{\partial }{\partial u}
 P_{0}[\eta _{(-u,v)}>t] \right) P[\tau_{D}(z)  >u]f^{*}(v)dudv\label{int-by-parts3}\\
 & \leq & \int_{0}^{\infty} \!\int_{0}^{\infty}  \left( \frac{\partial }{\partial u}
 P_{0}[\eta _{(-u,v)}>t] \right) P[\tau_{D_{*}}(0)  >v]f^{*}(u)dvdu \label{second-bm-isop}\\
 &=  &P_{0}[\tau_{D_{*}}(Z)>t]\label{int-by-parts4}
\end{eqnarray}
Equations (\ref{first-bm-isop}) and (\ref{second-bm-isop}) follow from the corresponding isoperimetric 
inequalities for Brownian motion and Lemma \ref{increasing-isop}.
Equations  (\ref{int-by-parts0}), (\ref{int-by-parts2}) and  (\ref{int-by-parts4}) follow by integration by parts, 
Lemma \ref{integration-by-parts-lemma}. Equation (\ref{int-by-parts3}) follows from Fubini-Tonelli theorem 
and integration by parts as in Lemma \ref{integration-by-parts-lemma}.

\end{proof}

We next give the proof of Theorem \ref{gen-isop-IBM}.
Let $X_{t}$ and $Y_{t}$
 be two  independent
one-dimensional Brownian motions, all started at $0$. BTBM is
defined to be $Z_{t}^{1}=x+X(|Y_{t}|)$. In $\RR{R}^{n}$, we
require $X$ to be independent one dimensional iterated Brownian
motions. 
Let $\tau_{D}(Z^{1})$ stand for the first exit time of BTBM from $D$. We have by the continuity of paths
\begin{equation}\label{probeqisop}
P_{z}[\tau _{D}(Z^{1}) > t]=P[\eta(-\tau _{D}(z) , \tau _{D}(z) ) >t].
\end{equation}

To prove isoperimetric-type inequalities we need the following
Lemma.
\begin{lemma}\label{isopheart}
$\frac{\partial}{\partial u} P_{0}[\eta _{(-u,u)}>t]\geq 0$ for all $u>0$.
\end{lemma}
\begin{proof}
Since $P_{0}[\eta _{(-u,u)}>t]\leq P_{0}[\eta _{(-v,v)}>t]$ for
$u\leq v$, the derivative is positive.

\end{proof}

\begin{proof}[Proof of Theorem \ref{gen-isop-IBM}]
 From equation (\ref{probeqisop})
\begin{eqnarray}
P_{z}[\tau _{D}(Z^{1}) > t]& = & P[\eta(-\tau _{D}(z) , \tau _{D}(z) ) >t]\nonumber \\
\ & \ = & \int_{0}^{\infty}P_{0}[\eta _{(-u,u)}>t]
f(u)du\label{genisopintegral}
\end{eqnarray}
By integration by parts, equation (\ref{genisopintegral}) equals
$$\int_{0}^{\infty}  \frac{\partial}{\partial u} P_{0} (\eta _{(-u,u)}>t)   P[\tau_{D}(z)>u]du$$
Now from Lemma \ref{isopheart}, $\frac{\partial }{\partial u} P_{0} (\eta
_{(-u,u)}>t)\geq 0$ for all $u>0$, which gives the desired
conclusion from the  corresponding isoperimetric inequalities for
the first exit time of the Brownian motion $X$ from $D_{*}$.

\end{proof}
We have a general result which follows from Theorems \ref{gen-isop-IBM} and \ref{gen-isop-burdzy-IBM}.
\begin{corollary}
Let $\xi $ and $T$ be  positive random variables such that $$\lim_{t\to\infty}P[\xi >t]=0=\lim_{t\to\infty}P[T >t]$$ 
and for each $t>0$
$$
 P[\xi >t]\leq P[T >t].
$$  Let $\xi_{1}, \xi_{2}$ be independent copies of $\xi$, $T_{1}, T_{2}$ independent copies of $T$ and let $Y$ be a one
dimensional Brownian motion independent of $\xi$ and $T$. Then for each $t>0$,
$$
P[\eta (-\xi_{1}, \xi_{2})>t]\leq P[\eta (-T_{1}, T_{2})>t]
$$
and
$$
P[\eta (-\xi, \xi)>t]\leq P[\eta (-T, T)>t].
$$
\end{corollary}

\textbf{Acknowledgments.} I would like to thank  Professor Rodrigo
 Ba\~{n}uelos, my academic advisor, for suggesting this problem to me and for  his guidance on this paper.

\end{document}